# $C^*$-ALGEBRAS OF HILBERT MODULE PRODUCT SYSTEMS

ILAN HIRSHBERG

ABSTRACT. We consider a class of $C^*$-algebras associated to one parameter continuous tensor product systems of Hilbert modules, which can be viewed as continuous counterparts of Pimsner's Toeplitz algebras. By exhibiting a homotopy of quasihomomorphisms, we prove that those algebras are $K$-contractible. One special case is closely related to the Rieffel-Wiener-Hopf extension of a crossed product by $\mathbb{R}$ considered by Rieffel and by Pimsner and Voiculescu, and can be used to produce a new proof of Connes' analogue of the Thom isomorphism and in particular of Bott periodicity. Another special case is closely related to Arveson's spectral $C^*$-algebras, and is used to settle Arveson's problem of computing their $K$-theory, extending earlier results of Zacharias to cover the general case.

Pimsner ([P]) introduced a class of $C^*$-algebras which is based on the left convolution operators acting on the 'Fock space' of a Hilbert module. Explicitly, for a given full right Hilbert $\mathcal{A}$-module $E$ along with a left action $\varphi : \mathcal{A} \to \mathcal{B}(E)$, one defines $\mathcal{E}_+ = \bigoplus_{n=0}^\infty E^{\otimes n}$. The elements of $E$ act on $\mathcal{E}_+$ by tensoring on the left, and then one lets $\mathcal{T}_E \subset \mathcal{B}(\mathcal{E}_+)$ be the $C^*$-algebra generated by those left tensoring operators. For our purposes, one can more conveniently think of $\mathcal{T}_E$ as the $C^*$-algebra generated by the images of $l^1(\mathcal{E}_+)$ in $\mathcal{B}(\mathcal{E}_+)$, where by $l^1(\mathcal{E}_+)$ we mean the sections $f$, with $f(n) \in E^{\otimes n}$, and $\sum_{n=0}^\infty \|f(n)\| < \infty$, and the action is given by convolution on the left. $\mathcal{E}_+$ can be thought of as a 'discrete product system', i.e. a fibration over $\mathbb{N}$ such that each fiber is a Hilbert module, along with a bilinear product structure which maps the fibers $E_n$ over $n$ and $E_m$ over $m$ to the fiber $E_{n+m}$ over $n+m$ and descends to an isomorphism $E_n \otimes_\mathcal{A} E_m \to E_{n+m}$. This presentation is unnecessarily complicated in this case, since the structure will be completely determined by the fiber over 1, but taking lead from Arveson's work on noncommutative dynamics ([A6], see example 2 below), it suggests a way to find a continuous analogue.

We shall consider an algebra $\mathcal{W}_E$ arising from the convolution action of $L^1(E)$ on a continuous analogue of the 'Fock space', $\int_{\mathbb{R}_+}^\oplus E_x dx$, where $E$ will be a product system modeled over $\mathbb{R}_+$ instead of $\mathbb{N}$. Our main result is the following.

**Theorem 1.** $KK(\mathcal{W}_E, \mathcal{W}_E) = 0$, i.e. $\mathcal{W}_E$ is $K$-contractible.

This can be viewed as the counterpart of the fact that $\mathcal{T}_E$ is $KK$-equivalent to $\mathcal{A}$ in the case of Pimsner's algebras. We shall use Cuntz's quasihomomorphism picture of $KK$-theory ([C1, C2], see also [B]).

**Definition 2.** Let $\mathcal{A}$ be a separable $C^*$-algebra. A *product system of $\mathcal{A}$-bimodules* $E$ is a measurable bundle of $\mathcal{A}$-bimodules over $\mathbb{R}_+$, along with a multiplication map $E \times E \to E$, which descends to an isomorphism $E_x \otimes_\mathcal{A} E_y \to E_{x+y}$ for all $x, y \in \mathbb{R}_+$ (where $E_x$ is the fiber over $x$), and is measurable in the sense that if $\xi$ is a measurable section and $e \in E_y$ then the sections $x \mapsto e\xi(x-y)$, $x \mapsto \xi(x-y)e$ (0 if $x < y$) are also measurable.

We refer the reader to the appendix for a discussion of measurable bundles of Hilbert bimodules.

The elements of $E$ act on $L^2(\mathbb{R}_+)$ on the left as adjoinable operators. Slightly abusing notation, we will use the same notation for $e \in E$ and $e$ thought of as an element of $\mathcal{B}(L^2(\mathbb{R}_+))$. Note that $\|e\|_{E_x} \geq \|e\|_{\mathcal{B}(\int_{\mathbb{R}_+}^\oplus E_x dx)}$. We can integrate those to get an action of $L^1(E)$ (the space of measurable sections $\xi$ such that $\int_{\mathbb{R}_+} \|\xi(x)\| dx < \infty$) on $\int_{\mathbb{R}_+}^\oplus E_x dx$.

**Definition 3.** For $f \in L^1(E)$ we define $W_f \in \mathcal{B}\left(\int_{\mathbb{R}_+}^\oplus E_x dx\right)$ by

$$W_f = \int_{\mathbb{R}_+} f(x) dx$$

(where here $f(x)$ is thought of as an element of $\mathcal{B}\left(\int_{\mathbb{R}_+}^\oplus E_x dx\right)$).

We denote by $\mathcal{W}_E$ the $C^*$-subalgebra of $\mathcal{B}\left(\int_{\mathbb{R}_+}^\oplus E_x dx\right)$ generated by

$$\{W_f \mid f \in L^1(E)\}$$

Examples:

(1) Let $\alpha$ be a one parameter group of automorphisms, or a semigroup of endomorphisms, acting on $\mathcal{A}$. Let $E = \mathcal{A} \times \mathbb{R}_+$, the right Hilbert $\mathcal{A}$ module structure being the usual one ($e \cdot a = ea$, $\langle a, b \rangle = a^*b$) and with the left action on the fiber $E_x$ being given by $a \cdot e_x = \alpha_x(a)e_x$. The product of $a \in E_x = \mathcal{A}$ by $b \in E_y = \mathcal{A}$ is $\alpha_y(a)b \in E_{x+y} = \mathcal{A}$.

$\mathcal{W}_E$ here is the algebra $\mathcal{T}_\alpha$ considered in [KS] (Definition 3.2). When $\alpha$ is an action of $\mathbb{R}$ by automorphisms, this algebra, in turn, is Morita equivalent to the Wiener-Hopf extension $\mathcal{W}_\alpha = C\mathcal{A} \times_{\tau \otimes \alpha} \mathbb{R}$ of $\mathcal{A} \times_\alpha \mathbb{R}$ by $\mathcal{A} \otimes \mathcal{K}$ considered by Rieffel and by Pimsner and Voiculescu in their 'Wiener-Hopf' proofs of Connes' analogue of the Thom isomorphism theorem ([R], [B] chapter 10; see [KS] for a proof of Morita equivalence).

Therefore by Theorem 1 the middle term in the sequence

$$0 \to \mathcal{A} \otimes \mathcal{K} \to \mathcal{W}_\alpha \to \mathcal{A} \times_\alpha \mathbb{R} \to 0$$

is $K$-contractible, and from that the Connes – Thom isomorphism

$$K_*(\mathcal{A}) \cong K_{*+1}(\mathcal{A} \times_\alpha \mathbb{R})$$

follows immediately by considering the long exact sequence in $K$-theory (where we can replace $K_*$ here by any half-exact stable homotopy functor from separable $C^*$-algebras to Abelian groups). Further specializing to the case of the trivial action, we recover Bott periodicity.

When $\alpha$ is a semigroup of endomorphisms, we obtain a different proof of the $K$-contractibility result from [KS] (Theorem 3.7).

(2) For $\mathcal{A} = \mathbb{C}$, there is an abundance of examples which are closely related to noncommutative dynamics ([A2, A3, A4, A5], see more recently [A6]). $\mathcal{W}_E$ here is Morita equivalent to Arveson's spectral $C^*$-algebra $C^*(E)$ ([A6] section 4.1). Theorem 1 shows then that $C^*(E)$ is $K$-contractible for any product system $E$, settling this problem of Arveson.

We note that Zacharias ([Z]) has proved that $C^*(E)$ is $K$-contractible whenever $E$ is a product system of Hilbert spaces which admits a non-zero measurable multiplicative section (product systems which are not of type $III$), using different techniques.

(3) Let $\mathcal{S}$ be a finite set. For $s \in \mathcal{S}$, we let $p_s \in C(\mathcal{S})$ be the function given by $p_s(t) = \delta_{s,t}$. Let $(\Omega, \mathfrak{B}, \mu)$ be a Lebesgue probability space. Let $X_x : \Omega \to \mathcal{S}$, $x \in [0, \infty)$ be a stationary Markov process, and let $\Phi_x : \Omega \to \Omega$ be a semigroup of measure preserving maps such that $X_x \circ \Phi_y = X_{x+y}$ for all $x, y \in [0, \infty)$. Denote by $\mathfrak{B}_{[x,y]}$ the $\sigma$-algebra generated by $\{X_r \mid x \leq r \leq y\}$. For $x \geq 0$, $s \in \mathcal{S}$, let $\Omega_x^s = X_x^{-1}\{s\}$. Note that $\bigcup_{s \in \mathcal{S}} \Omega_x^s = \Omega$ for all $x$. The Markov property implies that for any $s \in \mathcal{S}$, $x > 0$, if let $\mu_x^s$ be the restriction of $\frac{1}{\mu(\Omega_x^s)}\mu$ to $\Omega_x^s$, then $\{\Omega_x^s \cap Y \mid Y \in \mathfrak{B}_{[0,x]}\}$ and $\{\Omega_x^s \cap Y \mid Y \in \mathfrak{B}_{[x,\infty)}\}$ are independent with respect to the probability measure $\mu_x^s$. Therefore, for any $x < y < z$, pointwise multiplication induces an isomorphism

$$L^2\left(\Omega_y^s, \mathfrak{B}_{[x,y]}, \frac{1}{\mu(\Omega_y^s)}\mu\right) \otimes L^2\left(\Omega_y^s, \mathfrak{B}_{[y,z]}, \frac{1}{\mu(\Omega_y^s)}\mu\right) \cong L^2\left(\Omega_y^s, \mathfrak{B}_{[x,z]}, \frac{1}{\mu(\Omega_y^s)}\mu\right)$$

Let $E_{[x,y]} = L^2(\Omega, \mathfrak{B}_{[x,y]}, \mu)$ (as a vector space). We give $E_{[x,y]}$ the structure of a Hilbert bimodule over $C(\mathcal{S})$ as follows. For $s \in \mathcal{S}$, we let $p_s$ be the projection onto $L^2(\Omega_y^s, \mathfrak{B}_{[x,y]}, \mu)$. Note that $\sum_{s \in \mathcal{S}} p_s = 1$ under this action, which we think of as the right action of $C(\mathcal{S})$. The $C(\mathcal{S})$-valued inner product is given by

$$\langle \xi, \eta \rangle_{C(\mathcal{S})} = \sum_{s \in \mathcal{S}} \frac{1}{\mu(\Omega_y^s)} \langle \xi p_s, \eta p_s \rangle_{L^2} p_s$$

The left action $\varphi : C(\mathcal{S}) \to \mathcal{B}(E_{[x,y]})$ is given by letting $\varphi(p_s)$ be the projection onto $L^2(\Omega_x^s, \mathfrak{B}_{[x,y]}, \mu)$. Note that $\varphi(p_s)$ commutes with $p_t$ for all $s, t \in \mathcal{S}$, so the left action is given by module maps (which, in the case of Hilbert modules over $C(\mathcal{S})$, can easily be shown to be automatically adjoinable). As usual, we will suppress the $\varphi$, and just write the $p_s$ on the right or on the left to denote the left and right actions, when it causes no confusion.

A straightforward computation now shows that pointwise multiplication $E_{[x,y]} \times E_{[y,z]} \to E_{[x,z]}$ descends to an isomorphism $E_{[x,y]} \otimes_{C(\mathcal{S})} E_{[y,z]} \to E_{[x,z]}$.

Composition with $\Phi_x$ induces a one parameter semigroup of isometries $S_x$ of $L^2(\Omega, \mathfrak{B}, \mu)$. The $S_x$'s are furthermore module maps: when restricted to $E_{[y,z]}$, we have $S_x(\xi p_s) = S_x(\xi) p_s$ where the right action is taken in the appropriate sense in both cases.

Let $E_x = E_{[0,x]}$. We can view the bundle $E$ over $\mathbb{R}_+$ whose fiber over $x$ is $E_x$ as a sub-bundle of the trivial bundle $L^2(\Omega, \mathfrak{B}, \mu) \times \mathbb{R}_+$ over $\mathbb{R}_+$.

It clearly satisfies the conditions of being a measurable bundle of Hilbert bimodules. $E$ furthermore has the structure of a $C(\mathcal{S})$-product system. The idea is that starting with $E_x, E_y$, we use the shift $S_x$ to move $E_y$ to $E_{[x,x+y]}$, and then put together $E_x = E_{[0,x]}$ and the shifted version of $E_y$, $E_{[x,x+y]}$ to get $E_{x+y} = E_{[0,x+y]}$. Explicitly, we define a product $m : E_x \times E_y \to E_{x+y}$ by

$$m(\xi,\eta)(\omega) = \xi(\omega) \cdot (S_x\eta)(\omega) = \xi(\omega) \cdot \eta(\Phi_x(\omega)) \quad (\omega \in \Omega)$$

and it is now easy to see that $E$ is a $C(\mathcal{S})$-product system with this structure.

Were we to make this construction in the discrete case, we would obtain extensions of Cuntz-Krieger algebras (see [P], example 2), so $\mathcal{W}_E$ in this case can be thought of as a 'continuous time' analogue of a Cuntz-Krieger algebra.

## Proof of $K$-contractibility

The purpose of this section is to prove that $KK(\mathcal{W}_E, \mathcal{W}_E) = 0$. This will be done as follows. We define a family of maps $\pi_t : \mathcal{W}_E \to \mathcal{M}(\mathcal{K} \otimes \mathcal{W}_E)$, $t \in (0,\infty]$, such that the pairs $(\pi_\infty, \pi_t)$ are all equivalent quasihomomorphisms $\mathcal{W}_E \to \mathcal{K} \otimes \mathcal{W}_E$. We then show that their $KK$ class is equal both to the identity class $id_{\mathcal{W}_E}$ and to the 0 class. The former will be done by letting the parameter $t$ tend to 0, and the latter will be done by letting $t$ tend to $\infty$.

Let $\{S_x\}_{x \in \mathbb{R}_+}$ denote the semigroup of unilateral shifts on $L^2(\mathbb{R}_+)$.

**Definition 4.** We let $\pi_\infty : \mathcal{W}_E \to \mathcal{B}\left(L^2(\mathbb{R}_+) \otimes \int_{\mathbb{R}_+}^\oplus E_x dx\right)$ be the map given by

$$\pi_\infty(W_f) = \int_{\mathbb{R}_+} S_x \otimes f(x) dx$$

Here $L^2(\mathbb{R}_+) \otimes \int_{\mathbb{R}_+}^\oplus E_x dx$ denotes the exterior tensor product of those two modules (the first one thought of as a $\mathbb{C}$-module).

It is clear that $\pi_\infty$ as defined extends to a $*$-homomorphism from the $*$-subalgebra generated by the $W_f$'s. Showing that it extends to the closure $\mathcal{W}_E$ requires some argument. Note that for any $\lambda > 0$, the submodule $\int_{(\lambda,\infty)}^\oplus E_x dx$ is invariant for the action of the $W_f$. Furthermore, we have

$$\int_{(\lambda,\infty)}^\oplus E_x dx \cong \int_{\mathbb{R}_+}^\oplus E_x dx \otimes_\mathcal{A} E_\lambda$$

via the contraction map, by the definition of a product system. This isomorphism identifies the restriction of $W_f$ to $\int_{(\lambda,\infty)}^\oplus E_x dx$ with $W_f \otimes 1_{E_\lambda}$, and therefore the restriction extends to a homomorphism

$$\sigma_\lambda : \mathcal{W}_E \to \mathcal{B}\left(\int_{(\lambda,\infty)}^\oplus E_x dx\right)$$

(given by $\sigma_\lambda = id_{\mathcal{W}_E} \otimes 1_{E_\lambda}$).

For $\lambda \leq 0$, we denote for convenience $\int_{(\lambda,\infty)}^{\oplus} E_x dx = \int_{\mathbb{R}_+}^{\oplus} E_x dx$, with $\sigma_\lambda = id_{\mathcal{W}_E}$. Viewing $L^2(\mathbb{R}_+)$ as the direct integral of the trivial bundle over $\mathbb{R}_+$ with fiber $\mathbb{C}$, we can view $L^2(\mathbb{R}_+) \otimes \int_{\mathbb{R}_+}^{\oplus} E_x dx$ as a direct integral over the quarter plane $\mathbb{R}_+^2$, where the fiber over $(x,y)$ is $\mathbb{C}_x \otimes E_y$ (here $\mathbb{C}_x$ denotes a copy of $\mathbb{C}$).

For any $\lambda \in \mathbb{R}$, note that the restriction of this measurable bundle over the quarter plane to the ray $x - y = \lambda$ is isomorphic to the bundle $E|_{(\lambda,\infty) \cap \mathbb{R}_+}$. Extending the measurable bundle over the quarter plane to a measurable bundle over the half plane $y + x > 0$ by setting all fibers outside the quarter plane to be the 0 bimodule, we have a measurable bundle over $\mathbb{R} \times \mathbb{R}_+$. Thus, $\{\int_{(\lambda,\infty)}^{\oplus} E_x dx \mid \lambda \in \mathbb{R}\}$ can be given the structure of a measurable bundle over $\mathbb{R}$, and we have

$$L^2(\mathbb{R}_+) \otimes \int_{\mathbb{R}_+}^{\oplus} E_x dx \cong \int_{\mathbb{R}}^{\oplus} \int_{(\lambda,\infty)}^{\oplus} E_x dx d\lambda$$

Under this identification, we have $\pi_\infty(W_f) = \int_{\mathbb{R}}^{\oplus} \sigma_\lambda(\mathcal{W}_f) d\lambda$, and thus, $\pi_\infty$ extends to all of $\mathcal{W}_E$ by $\pi_\infty = \int_{\mathbb{R}}^{\oplus} \sigma_\lambda d\lambda$, as required (the family $\sigma_\lambda$ is measurable, since we can view the bundle over $\mathbb{R}$ with fiber $\int_{(\lambda,\infty)}^{\oplus} E_x dx$ as the interior tensor product of the trivial bundle over $\mathbb{R}$ with fiber $\int_{\mathbb{R}_+}^{\oplus} E_x dx$ by the bundle over $\mathbb{R}$ given by $E$ on $\mathbb{R}_+$ and by the trivial bundle with fiber $\mathcal{A}$ on $(\infty, 0]$; with this identification, we have $\sigma_\lambda \cong \sigma \otimes_\mathcal{A} 1$ throughout).

Note also that since $e_x \in \mathcal{M}(\mathcal{W}_E)$, we have that $\pi_\infty(\mathcal{W}_E) \subseteq \mathcal{M}(\mathcal{K} \otimes \mathcal{W}_E)$.

**Definition 5.** For $t > 0$ denote $h_t(x) = \sqrt{2t} e^{-tx}$.

Note that $h_t \in L^1(\mathbb{R}_+) \cap L^2(\mathbb{R}_+)$, $\|h_t\|_2 = 1$, $\|h_t\|_1 = \sqrt{2/t}$. Let $P_t$ denote the projection onto $h_t$ in $L^2(\mathbb{R}_+)$, and let $Q_t = 1 - P_t$. Since $S_x^* h_t = e^{-tx} h_t$, we have that for all $t$, $Q_t L^2(\mathbb{R}_+)$ is an invariant subspace for the semigroup $\{S_x\}$, and therefore $\{S_x Q_t\}_{x \in \mathbb{R}_+}$ is a semigroup of isometries of the space $Q_t L^2(\mathbb{R}_+)$ (which is unitarily equivalent to the unilateral shift semigroup).

**Definition 6.** For any $t > 0$, we let $\pi_t : \mathcal{W}_E \to \mathcal{B}\left(L^2(\mathbb{R}_+) \otimes \int_{\mathbb{R}_+}^{\oplus} E_x dx\right)$ be the map given by

$$\pi_t(W_f) = \int_{\mathbb{R}_+} S_x Q_t \otimes f(x) dx$$

$\pi_t$ is extends to all of $\mathcal{W}_E$ (by the same argument as that of $\pi_\infty$). As above, we also have $\pi_t(\mathcal{W}_E) \subseteq \mathcal{M}(\mathcal{K} \otimes \mathcal{W}_E)$.

**Definition 7.** Let $\varphi_t : \mathcal{A} \to \mathcal{K} \otimes \mathcal{W}_E \subseteq \mathcal{B}\left(L^2(\mathbb{R}_+) \otimes \int_{\mathbb{R}_+}^{\oplus} E_x dx\right)$ denote the homomorphism $\varphi_t(T) = P_t \otimes T$.

$\varphi_t$ represents, for any $t$, the identity class $id_{\mathcal{W}_E}$ in $KK(\mathcal{W}_E, \mathcal{W}_E)$. Note that $\pi_t, \varphi_t$ have orthogonal supports, so $\pi_t + \varphi_t$ is also a homomorphism.

**Claim 8.** $\|S_x P_t - S_{x+\delta} P_t\| = \sqrt{2 - 2e^{-\delta t}}$ for any $\delta > 0$.

*Proof.*
$$\|S_x P_t - S_{x+\delta} P_t\| = \|P_t - S_\delta P_t\| = \|(S_\delta - 1)h_t\|_2$$
and
$$\|(S_\delta - 1)h_t\|_2^2 = 2t \int_0^\delta e^{-2ty} dy + 2t(e^{t\delta} - 1)^2 \int_\delta^\infty e^{-2ty} dy =$$
$$= 1 - e^{-2t\delta} + (e^{t\delta} - 1)^2 e^{-2t\delta} = 2 - 2e^{-t\delta}$$
$\square$

**Claim 9.** *For any $T \in \mathcal{W}$, $\pi_\infty(T) - \pi_t(T) \in \mathcal{K} \otimes \mathcal{W}_E$.*

*Proof.* It suffices to verify it for $T$ of the form $W_f$, $f \in L^1(E)$. Fix $t$. Fix $\epsilon > 0$. Choose $\delta > 0$ such that $\|f\|_1 \sqrt{2 - 2e^{-\delta t}} < \epsilon$. It suffices to show that we can approximate our element arbitrarily by elements from $\mathcal{K} \otimes \mathcal{W}_E$. So,
$$\pi_\infty(W_f) - \pi_t(W_f) = \int_{\mathbb{R}_+} (S_x \otimes f(x) - S_x Q_t \otimes f(x)) dx =$$
$$= \int_{\mathbb{R}_+} (S_x P_t \otimes f(x)) dx = \sum_{n=0}^\infty \int_{n\delta}^{(n+1)\delta} (S_x P_t \otimes f(x)) dx =$$
$$= \sum_{n=0}^\infty \int_{n\delta}^{(n+1)\delta} (S_{n\delta} P_t \otimes f(x)) dx - R$$
Where
$$R = \sum_{n=0}^\infty \int_{n\delta}^{(n+1)\delta} ((S_{n\delta} - S_x) P_t \otimes f(x)) dx$$
Note that
$$\|R\| \leq \sum_{n=0}^\infty \int_{n\delta}^{(n+1)\delta} \|((S_{n\delta} - S_x) P_t \otimes f(x)\| dx \leq \sum_{n=0}^\infty \int_{n\delta}^{(n+1)\delta} \|f(x)\| \epsilon / \|f\|_1 dx = \epsilon$$
and
$$\sum_{n=0}^\infty \int_{n\delta}^{(n+1)\delta} (S_{n\delta} P_t \otimes f(x)) dx = \sum_{n=0}^\infty S_{n\delta} P_t \otimes W_{f \chi_{[n\delta,(n+1)\delta]}} \in \mathcal{K} \otimes \mathcal{W}_E$$
(the last sum converges in norm). $\square$

**Corollary 10.** *The pairs $\Phi_t = (\pi_\infty, \pi_t)$ and $\Psi_t = (\pi_\infty, \pi_t + \varphi_t)$ are quasihomomorphisms $\mathcal{W}_E \to \mathcal{K} \otimes \mathcal{W}_E$. Note that $[\Psi_t] = [\Phi_t] - [id_{\mathcal{W}_E}]$. $t \mapsto \pi_\infty(T) - \pi_t(T)$ is norm-continuous for any fixed $T \in \mathcal{W}_E$, $t \in (0, \infty)$. Therefore $\Phi_t, \Psi_t$ are homotopies of quasihomomorphisms.*

**Lemma 11.** *$[\Psi_t] = 0$ in $KK(\mathcal{W}_E, \mathcal{W}_E)$, i.e. $[\Phi_t] = [id_{\mathcal{W}_E}]$ for all $t \in (0, \infty)$.*

*Proof.* We will show that for any $T \in \mathcal{W}_E$, $\|\pi_\infty(T) - \pi_t(T) - \varphi_t(T)\| \to 0$ as $t \to 0$. It suffices to check this for $T$ in a set of generators, and so it suffices to check this for all $T = W_f$ such that $f$ has compact support. Suppose $supp(f) \subseteq [0, M]$.
$$\|\pi_\infty(W_f) - \pi_t(W_f) - \varphi_t(W_f)\| = \left\|\int_{\mathbb{R}_+} ((S_x - 1) P_t \otimes f(x)) dx\right\|$$

From Claim 8 we know that $\|(S_x - 1)P_t\| = \sqrt{2 - 2e^{-tx}}$. The above integral then is:

$$\left\| \int_0^M ((S_x - 1)P_t \otimes f(x))dx \right\| \leq \int_0^M \|(S_x - 1)P_t \otimes f(x)\|dx \leq$$

$$\leq \int_0^M \|f(x)\|\sqrt{2 - 2e^{-tx}}dx \leq \|f\|_1 \sqrt{2 - 2e^{-tM}} \to 0$$

as $t \to 0$. $\square$

In order to prove the other part, we first need the following inequality.

**Lemma 12.** *Let $H$ be a Hilbert space, $\mathcal{E}$ be a Hilbert $\mathcal{A}$-module, $(\Omega, \mathfrak{B}, \mu)$ a measure space and let $v : \Omega \to H$, $w : \Omega \to \mathcal{E}$ be measurable functions with $\|w(\omega)\|_\mathcal{E} \leq 1$ a.e., and such that $\langle v(\omega), v(\eta) \rangle \geq 0$ a.e. If $\|v(\omega)\|_H \in L^1(\Omega, \mu)$ then*

$$\left\| \int_\Omega v(\omega) \otimes w(\omega) d\mu \right\|_{H \otimes \mathcal{E}} \leq \left\| \int_\Omega v(\omega) d\mu \right\|_H$$

*Proof.*

$$\left\| \int_\Omega v(\omega) \otimes w(\omega) d\mu_\omega \right\|_{H \otimes \mathcal{E}}^2 = \left\| \int_{\Omega \times \Omega} \langle v(\omega) \otimes w(\omega), v(\eta) \otimes w(\eta) \rangle d(\mu \times \mu)_{\omega,\eta} \right\|_\mathcal{A} =$$

$$= \left\| \int_{\Omega \times \Omega} \langle v(\omega), v(\eta) \rangle \langle w(\omega), w(\eta) \rangle d(\mu \times \mu)_{\omega,\eta} \right\|_\mathcal{A} \leq$$

$$\leq \int_{\Omega \times \Omega} \langle v(\omega), v(\eta) \rangle \|\langle w(\omega), w(\eta) \rangle\|_\mathcal{A} d(\mu \times \mu)_{\omega,\eta} \leq$$

$$\leq \int_{\Omega \times \Omega} \langle v(\omega), v(\eta) \rangle d(\mu \times \mu)_{\omega,\eta} = \left\| \int_\Omega v(\omega) d\mu \right\|_H^2$$

$\square$

**Lemma 13.** $[\Phi_t] = 0$ *in* $KK(\mathcal{W}_E, \mathcal{W}_E)$.

*Proof.* We will show that for any $T \in \mathcal{W}_E$, $\|\pi_\infty(T) - \pi_t(T)\| \to 0$ as $t \to \infty$. It suffices to check this for $T = W_f$, $f \in L^1(E)$. We can furthermore assume that without loss of generality that $\int_{\mathbb{R}_+} \|f(x)\|^2 dx < \infty$.

$$\|\pi_\infty(W_f) - \pi_t(W_f)\| =$$

$$= \sup \left\{ \left\| \int_{\mathbb{R}_+} (S_x P_t \otimes f(x))\xi dx \right\|_{L^2(\mathbb{R}_+) \otimes \int_{\mathbb{R}_+}^\oplus E_x dx} \,\middle|\, \xi \in L^2(\mathbb{R}_+) \otimes \int_{\mathbb{R}_+}^\oplus E_x dx, \|\xi\| = 1 \right\} =$$

$$= \sup \left\{ \left\| \int_{\mathbb{R}_+} (S_x P_t \otimes f(x))(h_t \otimes \eta) dx \right\|_{L^2(\mathbb{R}_+) \otimes \int_{\mathbb{R}_+}^\oplus E_x dx} \,\middle|\, \eta \in \int_{\mathbb{R}_+}^\oplus E_x dx, \|\eta\| = 1 \right\}$$

$$= \sup \left\{ \left\| \int_{\mathbb{R}_+} \|f(x)\| S_x h_t \otimes \frac{1}{\|f(x)\|} f(x) \eta dx \right\|_{L^2(\mathbb{R}_+) \otimes \int_{\mathbb{R}_+}^\oplus E_x dx} \,\middle|\, \eta \in \int_{\mathbb{R}_+}^\oplus E_x dx, \|\eta\| = 1 \right\}$$

where if $f(x) = 0$ then $\frac{1}{\|f(x)\|}f(x)\eta$ is taken to mean 0.

The families of vectors $v(x) = \|f(x)\|S_x h_t$, $w(x) = \frac{1}{\|f(x)\|}f(x)\eta$ satisfy the conditions of Lemma 12. So, the last expression is bounded by

$$\left\|\int_{\mathbb{R}_+} \|f(x)\|S_x h_t dx\right\|_{L^2(\mathbb{R}_+)}$$

Denote $g(x) = \|f(x)\|$. The last expression is equal to

$$\|g * h_t\|_{L^2(\mathbb{R}_+)}$$

Now, since $g \in L^1(\mathbb{R}_+) \cap L^2(\mathbb{R}_+)$, we have

$$\|g * h_t\|_{L^2} \leq \|g\|_{L^2}\|h_t\|_{L^1} = \sqrt{2}\|g\|_{L^2}/\sqrt{t} \to 0$$

as $t \to \infty$ □

This concludes the proof of Theorem 1.

## Appendix: direct integrals of Hilbert modules

In this appendix we discuss direct integrals of Hilbert bimodules. This is mostly a straightforward generalization of the well established theory for Hilbert spaces (see [D]), and we thus omit or sketch most of the proofs.

## Measurable bundles of bimodules and direct integrals

Let $(\Omega, \mathfrak{B}, \mu)$ be a $\sigma$-finite measure space.

**Definition A.1.** Let $\mathcal{A}$ be a separable $C^*$-algebra. A *measurable bundle of Hilbert $\mathcal{A}$-bimodules over $\Omega$*, $E$, is a collection $\{E_x \mid x \in \Omega\}$ of right Hilbert $\mathcal{A}$-modules with left actions via adjoinable operators, along with a distinguished vector subspace $\Gamma$ of $\Pi_{x \in \Omega} E_x$ (called the set of *measurable sections*) such that
  (1) For any $\xi \in \Gamma$, $a \in \mathcal{A}$, the functions $x \mapsto \langle \xi(x), \xi(x) \rangle$, $x \mapsto \langle a\xi(x), \xi(x) \rangle$ are measurable (as functions $\Omega \mapsto \mathcal{A}$).
  (2) If $\eta \in \Pi_{x \in \Omega} E_x$ satisfies that $x \mapsto \langle \xi(x), \eta(x) \rangle$ is measurable for all $\xi \in \Gamma$ then $\eta \in \Gamma$.
  (3) There exists a countable subset $\xi_1, \xi_2, ...$ of $\Gamma$ such that for all $x \in \Omega$, $\xi_1(x), \xi_2(x), ...$ are dense in $E_x$.

**Remark A.2.** We need to be precise about the measurable structure we use on $\mathcal{A}$. $\mathcal{A}$ is a complete separable metric space with its norm, and we use that topology to generate a measurable $\sigma$-algebra of Borel sets (hence it is a standard Borel space). Note that if we use the weak topology instead of the norm topology, we get the same Borel structure (see [A1] Theorem 3.3.5).

**Remark A.3.** The following properties all follow easily from Definition A.1.
  (1) Suppose that $\xi \in \Gamma$, and $f : \Omega \to \mathbb{C}$ is measurable, then the section $f \cdot \xi$ given by $f \cdot \xi(x) = f(x)\xi(x)$ is measurable.
  (2) If $\xi \in \Gamma$, $a \in \mathcal{A}$ then $x \mapsto \xi(x)a$ is also in $\Gamma$.

(3) If $\xi, \eta \in \Gamma$ then a simple polarization argument shows that $x \mapsto \langle \xi(x), \eta(x) \rangle$ is measurable.
(4) If $\xi \in \Gamma$, $a \in \mathcal{A}$ then $x \mapsto a\xi(x)$ is also in $\Gamma$.

**Lemma A.4.** *Suppose we have a collection $\{E_x\}_{x \in \Omega}$ as above, and suppose $\Gamma^0$ is a countable collection of sections such that*

(1) *If $\xi, \eta \in \Gamma^0$ then $x \mapsto \langle \xi(x), \eta(x) \rangle$ is measurable.*
(2) *If $\xi, \eta \in \Gamma^0$, $a \in \mathcal{A}$ then $x \mapsto \langle a\xi(x), \eta(x) \rangle$ is measurable.*
(3) *For all $x \in \Omega$, $\{\xi(x) \mid \xi \in \Gamma^0\}$ is dense in $E_x$.*

*then the set*

$$\Gamma = \{\xi \in \Pi_{x \in \Omega} E_x \mid x \mapsto \langle \xi(x) + a\xi(x), \eta(x) \rangle \text{ are measurable } \forall a \in \mathcal{A}, \, \eta \in \Gamma^0\}$$

*satisfies the condition of being the collection of measurable sections from Definition A.1. $\Gamma$ will be referred to as the* closure *of $\Gamma^0$.*

We now make precise what we mean by $\mathcal{A}$-valued integrals. We use weak integrals here. Let $f : \Omega \to \mathcal{A}$ be a measurable function. We say that $f$ is *integrable* if $\rho \mapsto \int_\Omega \rho(f(x)) d\mu$ is a well defined weak$^*$-continuous linear functional on $\mathcal{A}^*$ (the dual space of $\mathcal{A}$). This implies that there is a unique element $a \in \mathcal{A}$ which satisfies $\rho(a) = \int_\Omega \rho(f(x)) d\mu$ for all $\rho \in \mathcal{A}^*$, and we write $a = \int_\Omega f(x) d\mu$. Note that if $\int_\Omega \|f(x)\| dx < \infty$ then by the dominated convergence theorem, $f$ is integrable.

**Definition A.5.** We denote by $\Gamma^2 = \Gamma^2(E)$ the set of (a.e. equivalence classes of) measurable sections $\xi$ of $E$, such that $\int_\Omega \|f(x)\|^2 dx < \infty$. It is easy to see that $\Gamma^2$ is a pre-Hilbert $\mathcal{A}$-module, and we denote by $\int_\Omega^\oplus E_x d\mu$ its completion.

Note that if $\xi \in \Gamma^2$ and $a \in \mathcal{A}$ then $x \mapsto a\xi(x)$ is also in $\Gamma^2$, and is clearly adjoinable (with adjoint given by the action of $a^*$). Therefore $\int_\Omega^\oplus E_x d\mu$ has a left action of $\mathcal{A}$.

### Double direct integrals

Suppose $E$ is a measurable bundle of Hilbert $\mathcal{A}$-bimodules over $\Omega = \Omega_1 \times \Omega_2$. For any $x \in \Omega_1$ we have that $E|_{\{x\} \times \Omega_2}$ is naturally a measurable bundles of Hilbert $\mathcal{A}$-bimodules over $\Omega_2$. Let $F_x = \int_{\Omega_2}^\oplus E_{x,y} d\mu_y$, and let $F = \bigcup_{x \in \Omega_1} F_x$. We wish to define a collection $\Gamma$ of measurable sections which will make $F$ a measurable bundle of Hilbert $\mathcal{A}$-bimodules over $\Omega_1$.

To this end, we need some restrictions on the measure spaces involved. *From now on, we shall assume that each of our measure spaces admits a countable cover $\mathfrak{C}$ of measurable subsets of finite measure such that for any measurable set of finite measure $S$ and any $\epsilon > 0$ there exists a set $S' \in \mathfrak{C}$ such that the measure of the symmetric difference $(S - S') \cup (S' - S')$ is less than $\epsilon$.* This holds, for example, for any open subset of $\mathbb{R}$ with Lebesgue measure, where we can pick our collection $\mathfrak{C}$ to be the set of all finite unions of bounded intervals with rational endpoints. This is the only example we need in this paper.

Let $\xi_1, \xi_2, ...$ be sections of $E$ as in the definition of a measurable bundle (part 3). We may assume without loss of generality that the $\xi_i$'s are all bounded and vanish outside a set of finite measure, since otherwise we choose a countable disjoint cover $\Theta_1, \Theta_2, ...$ of $\Omega$ by measurable sets of finite measure, and replace the given collections by the collection

$$\chi_{\Theta_k} \cdot \chi_{\{x \mid m \leq \|\xi_n(x)\| < m+1\}} \cdot \xi_n \quad k, n, m \in \mathbb{N}$$

Note that in particular, this implies that $\xi_n \in \Gamma^2$ for all $n$.

Now suppose $\Omega = \Omega_1 \times \Omega_2$, and $\mathfrak{C}_1, \mathfrak{C}_2$ countable covers for $\Omega_1, \Omega_2$ as above. We take all finite unions of elements from $\{S_1 \times S_2 \mid S_1 \in \mathfrak{C}_1, S_2 \in \mathfrak{C}_2\}$ as our countable cover $\mathfrak{C}$ for $\Omega$, and let $\tilde{\Gamma}$ denote the all finite sums of sections of the form $\chi_S \xi_n$, $S \in \mathfrak{C}$ and $\xi_n$ as above. Note that this collection is countable and also satisfies the conditions of the definition of a measurable bundle, so we could take this to be our $\xi_i$'s. For any $x \in \Omega_1$, we denote by $\tilde{\Gamma}_x$ the collection of all restrictions of elements of $\tilde{\Gamma}$ to $x \times \Omega_2$, thought of as a subset of $F_x$. One checks that $\tilde{\Gamma}_x$ is dense in $F_x$.

Any section $\xi \in \tilde{\Gamma}$ gives rise to a section of $F$, namely $x \mapsto \xi|_{\{x\} \times \Omega_2}$. Denote this section by $\bar{\xi}$, and the set of all sections arising in this manner by $\overline{\Gamma^0}$. We therefore can view $F$ as a measurable bundle of Hilbert $\mathcal{A}$-bimodules over $\Omega_1$, with the measurable sections given by the closure of $\overline{\Gamma^0}$. The following 'Fubini theorem' now follows immediately.

**Theorem A.6** ('Fubini's theorem' for direct integrals). *With the notation above,* $\int_{\Omega_1 \times \Omega_2}^{\oplus} E_{x,y} d(\mu_1 \times \mu_2) \cong \int_{\Omega_1}^{\oplus} F_x d\mu_1$

## Measurable bundles of operators

Let $(\Omega_1, \mu_1), (\Omega_2, \mu_2)$ be two $\sigma$-finite measure spaces, and let $E, F$ be measurable bundles of Hilbert $\mathcal{A}$-bimodules over $\Omega_1, \Omega_2$ respectively. Suppose $\theta : \Omega_1 \to \Omega_2$ is a 1-1 measure preserving measurable embedding of $\Omega_1$ into $\Omega_2$ (so for any measurable $\Xi \in \Omega_1$, we require that $\theta(\Xi)$ is measurable as well, and $\mu_2(\theta(\Xi)) = \mu_1(\Xi)$). A family of adjoinable operators $T_x : E_x \to F_{\theta(x)}$, $x \in \Omega_1$, denoted $T : E \mapsto F$, is said to be a *measurable family* if for any measurable section $\xi$ of $E$, the section $y \mapsto T_{\theta^{-1}(y)} \xi(\theta^{-1}(y))$ (understood to be 0 when $y \notin \theta(\Omega_1)$) is a measurable section of $F$.

We observe that if $x \mapsto \|T_x\|$ is bounded then this measurable family gives rise in a natural way to an adjoinable map $\int_{\Omega_1}^{\oplus} T d\mu_1 : \int_{\Omega_1}^{\oplus} E_x d\mu_1 \to \int_{\Omega_2}^{\oplus} F_y d\mu_2$, whose norm is bounded by $\sup\{\|T_x\|\}_{x \in \Omega_1}$.

## Tensor products

**Exterior tensor product.** Let $E, \Gamma_1$ be a measurable bundle of Hilbert $\mathcal{A}$-modules over $\Omega_1$, and let $F, \Gamma_2$ be a measurable bundle of Hilbert $\mathcal{B}$-modules over $\Omega_2$. We denote by $\mathcal{A} \otimes \mathcal{B}$ the maximal tensor product. $E \otimes F$ will naturally have the structure of a measurable bundle of Hilbert $\mathcal{A} \otimes \mathcal{B}$-modules over $\Omega_1 \times \Omega_2$ where we take $(E \otimes F)_{x,y} = E_x \otimes F_y$ (the exterior tensor product), and $\Gamma$ to be the closure of all sections of the form $(x, y) \mapsto \xi(x) \otimes \eta(y)$, $\xi \in \Gamma_1$, $\eta \in \Gamma_2$.

Suppose $E'$ is another measurable bundle of Hilbert $\mathcal{A}$-modules over $\Omega_1$, and $F'$ is another measurable bundle of Hilbert $\mathcal{B}$-modules over $\Omega_2$. If $T : E \to E'$ and $S: F \to F'$ are measurable families of adjoinable operators then we have a natural family of adjoinable operators $T \otimes S: E \otimes F \to E' \otimes F'$. If $T, S$ are bounded then so is $T \otimes S$.

**Interior tensor product.** Let $(E, \Gamma_1), (F, \Gamma_2)$ be two measurable bundles of Hilbert $\mathcal{A}$-modules over $\Omega$. We can form a measurable bundle $E \otimes_{\mathcal{A}} F$ over $\Omega$ given by $(E \otimes_{\mathcal{A}} F)_x = E_x \otimes_{\mathcal{A}} F_x$ and $\Gamma$ being the closure of all sections of the form $x \mapsto \xi(x) \otimes_{\mathcal{A}} \eta(x)$, $\xi \in \Gamma_1$, $\eta \in \Gamma_2$.

Suppose $E'$ is another measurable bundle of Hilbert $\mathcal{A}$-modules over $\Omega$, and $T: E \to E'$ is a measurable family of adjoinable operators then we have a natural family of adjoinable operators $T \otimes_{\mathcal{A}} 1: E \otimes_{\mathcal{A}} F \to E' \otimes_{\mathcal{A}} F$.

Note that $1 \otimes_{\mathcal{A}} T: F \otimes_{\mathcal{A}} E \to F \otimes_{\mathcal{A}} E'$ will be well defined if a.e. $x$, $T_x$ intertwines the left actions of $\mathcal{A}$ on $E_x$ and $E'_x$.

## References

[A1] Arveson, W.B., *An Invitation to $C^*$-Algebras*, 2nd printing, Graduate Texts in Mathematics **39**, Springer-Verlag, 1998.

[A2] Arveson, W.B., Continuous analogues of Fock space, *Mem. Amer. Math. Soc.* **80** (1989), no. 409.

[A3] Arveson, W.B., Continuous analogues of Fock space. II. The spectral $C^*$-algebra, *J. Funct. Anal.* **90** (1990), no. 1, 138–205.

[A4] Arveson, W.B., Continuous analogues of Fock space. III. Singular states. *J. Operator Theory* **22** (1989), no. 1, 165–205.

[A5] Arveson, W.B., Continuous analogues of Fock space. IV. Essential states. *Acta Math.* **164** (1990), no. 3-4, 265–300.

[A6] Arveson, W.B., *Noncommutative Dynamics and E-semigroups*, to appear in Springer Monographs in Mathematics, Springer-Verlag.

[B] Blackadar, B., *K-Theory for Operator Algebras*, 2nd edition, MSRI publications **5**, Cambridge University Press, 1998.

[C1] Cuntz, J., Generalized homomorphisms between $C^*$-algebras and $KK$-theory. *Dynamics and processes (Bielefeld, 1981)*, 31–45, Lecture Notes in Math. **1031**, Springer, Berlin, 1983.

[C2] Cuntz, J., $K$-theory and $C^*$-algebras. *Algebraic K-theory, number theory, geometry and analysis (Bielefeld, 1982)*, 55–79, Lecture Notes in Math., **1046**, Springer, Berlin, 1984.

[D] Dixmier, J., *Von-Neumann Algebras*, North Holland, 1981.

[P] Pimsner, M.V., A class of $C^*$-algebras generalizing both Cuntz-Krieger algebras and crossed products by $\mathbb{Z}$. *Free probability theory (Waterloo, ON, 1995)*, 189–212, Fields Inst. Commun., **12**, Amer. Math. Soc., Providence, RI, 1997.

[R] Rieffel, M., Connes' analogue for crossed products of the Thom isomorphism, *Operator algebras and K-theory (San Francisco, Calif., 1981)*, pp. 143–154, Contemp. Math., **10**, Amer. Math. Soc., Providence, R.I., 1982.

[KS] Khoshkam, M. and Skandalis, G., Toeplitz algebras associated with endomorphisms and Pimsner-Voiculescu exact sequences, *Pacific J. Math.*, **181** (1997) no. 2, 315–331.

[Z] Zacharias, J., Continuous tensor products and Arveson's spectral $C^*$-algebras, *Mem. Amer. Math. Soc.* **143** (2000), no. 680.

*E-mail address*: ilan@math.berkeley.edu

DEPARTMENT OF MATHEMATICS, UNIVERSITY OF CALIFORNIA, BERKELEY CA 94720, USA